\documentclass{article}
\usepackage[margin=1.2in]{geometry}

\usepackage{amsthm, amsmath, amssymb, inputenc,  comment, tikz, array, longtable, enumerate, graphicx, float, hyperref}
\usepackage{listings}
\usepackage{xcolor, algorithm, algpseudocode}
\usepackage{thmtools}
\usepackage{thm-restate}

\lstdefinestyle{mystyle}{
    backgroundcolor=\color{lightgray!20}, 
    commentstyle=\color{green!50!black},
    keywordstyle=\color{blue},
    stringstyle=\color{red!70!black},
    basicstyle=\ttfamily\small, 
    breaklines=true,             
    frame=single,                
    numbers=left,                
    numberstyle=\tiny\color{gray},
    stepnumber=1,
    tabsize=4,
    showstringspaces=false
}

\usepackage{amsthm, amsmath, amssymb, inputenc, authblk, comment, tikz, array, longtable, enumerate, comment, xcolor}
\usepackage{xcolor}
\usepackage{graphicx}
\usepackage{float}

\newtheorem{theorem}{Theorem}[section]
\newtheorem{lemma}[theorem]{Lemma}

\newtheorem{conjecture}[theorem]{Conjecture}
\newtheorem{question}[theorem]{Question}

\theoremstyle{definition}

\newtheorem{claim}{Claim}

\tikzset{every node/.style={inner sep=1.5pt,circle,draw,fill}}

\tikzstyle{every path}=[thick]

\newdimen\R
\title{Weighted Max Cut in Subcubic Graphs}

\author{
Will Hausmann\thanks{Partially supported by the U.S. Department of Education Graduate Assistance in Areas of National Need (GAANN) Award \# P200A240169.}\\
School of Mathematics, Georgia Institute of Technology, USA
}

\date{\today}

\begin{document}

\maketitle

\begin{abstract}
For a graph $G$ with nonnegative edge weights, let $\operatorname{mc}(G,w)$ denote the maximum total weight of a cut and let $w(G)$ denote the total edge weight.
We prove that every edge-weighted triangle-free subcubic graph satisfies $\operatorname{mc}(G,w)\geq \frac{15}{19}w(G)$, improving the bound of $\frac{8}{11}$ due to Gutin and Yeo.
For subcubic graphs of odd girth at least seven and nine, we strengthen this bound to $\frac{47}{59}w(G)$ and $\frac45w(G)$, respectively, establishing the conjecture of Gutin and Yeo for graphs of odd girth at least nine.
Our proofs combine structural properties of vertex partitions with randomized cut constructions and convex combinations of the resulting bounds.

\end{abstract}

\section{Introduction}

Let $G=(V,E,w)$ be a finite simple edge-weighted graph with weight function
$w:E\to[0,\infty)$. The \textit{girth} of $G$ is the length of its shortest cycle. The \textit{odd girth} of $G$ is the length of its shortest odd cycle. For a partition $(L,R)$ of $V(G)$, the
\textit{cut} $(L,R)$ is the spanning bipartite subgraph consisting
of the edges between $L$ and $R$. For disjoint vertex sets $A,B\subseteq V(G)$, let $E(A,B)$
denote the set of edges with one endpoint in $A$ and the other
in $B$, and let $G[A,B]$ denote the spanning bipartite subgraph
with edge set $E(A,B)$.
Given a cut $(L,R)$, an edge is \textit{internal} if it belongs
to $E(G[L])\cup E(G[R])$ and \textit{external} or
\textit{crossing} if it belongs to $E(L,R)$. The Max Cut problem asks for a
cut of maximum total weight, which we denote by
$\operatorname{mc}(G,w)$. In the unweighted case, every edge has
weight one, and we write $\operatorname{mc}(G)$ for the maximum
number of edges in a cut.

A graph is \textit{cubic} if every vertex has degree exactly three, and \textit{subcubic} if every vertex has degree at most three.
Computing a maximum cut is NP-hard even for triangle-free cubic
graphs~\cite{Yannakakis1978}. Moreover, for every $\varepsilon>0$,
a polynomial-time approximation algorithm for Max Cut on cubic
graphs with approximation ratio $\frac{331}{332}+\varepsilon$
would imply $\mathrm{P}=\mathrm{NP}$~\cite{BermanKarpinski1999}.
On the positive side, Halperin, Livnat, and
Zwick~\cite{HalperinLivnatZwick2004} gave a
$0.9326$-approximation algorithm for unweighted Max Cut on
subcubic graphs, using semidefinite programming and local
improvement, and a combinatorial $\frac{22}{27}$-approximation
algorithm for cubic graphs. These approximation ratios compare
the cut returned by an algorithm with an optimal cut; our focus
is on lower bounds relative to the total edge weight.

Subcubic graphs combine computational difficulty with structural
restrictions that yield strong lower bounds. Every subcubic graph
$G$ satisfies
\(
\operatorname{mc}(G)\geq\frac23|E(G)|,
\)
with equality for $K_3$ and $K_4$. Hopkins and
Staton~\cite{Hopkins1982} proved that every triangle-free cubic
graph has a cut containing at least $\frac45$ of its edges.
Bondy and Locke~\cite{BondyLocke1986} extended this bound to
triangle-free subcubic graphs. Xu and Yu~\cite{XuYu2008}
characterized the seven connected triangle-free subcubic graphs
for which equality holds. Zhu~\cite{Zhu2009} subsequently proved
that every 2-connected triangle-free subcubic graph outside the
exceptional family of Xu and Yu satisfies
\(
\operatorname{mc}(G)\geq\frac{17}{21}|E(G)|.
\)

For weighted graphs, Poljak and Turz\'ik~\cite{PoljakTurzik1986}
proved the general lower bound
\(
\operatorname{mc}(G,w)
\geq
\frac12w(G)+\frac14w(T_{\min}),
\)
where $G$ is connected and $T_{\min}$ is a minimum-weight spanning
tree of $G$. Gutin and Yeo~\cite{Gutin2023} developed a structural
approach that strengthens the inequality of Poljak and Turz\'ik
and yields the following bound for triangle-free subcubic graphs.

\begin{theorem}[Gutin and Yeo~\cite{Gutin2023}]
\label{thm:gutin-yeo-eight-elevenths}
Let $G=(V,E,w)$ be an edge-weighted triangle-free subcubic graph.
Then
\[
\operatorname{mc}(G,w)\geq\frac{8}{11}w(G).
\]
\end{theorem}

They also conjectured that the unweighted $\frac45$ bound extends
to arbitrary nonnegative edge weights.

\begin{conjecture}[Gutin and Yeo~\cite{Gutin2023}]
\label{conj: maxcut conj}
Let $G=(V,E,w)$ be an edge-weighted triangle-free subcubic graph.
Then
\[
\operatorname{mc}(G,w)\geq\frac45w(G).
\]
\end{conjecture}

A stronger weighted bound is implicit in a result of
Kardo\v{s}, Kr\'al', and Volec, stated in terms of probability
distributions on cuts rather than weighted Max Cut.

\begin{theorem}[Kardo\v{s}, Kr\'al', and Volec~\cite{KardosKralVolec2012}]
\label{thm:KKV-large-girth}
There exists an absolute constant $g_0$ such that every cubic
graph $G$ of girth at least $g_0$ admits a probability distribution
on cuts under which each edge crosses with probability at least
$0.88672$.
\end{theorem}

Their proof shows that one may take $g_0=637{,}789$.
For a fixed graph $G$,
this distributional statement is equivalent to the weighted
Max Cut bound
\(
\operatorname{mc}(G,w)\geq 0.88672\,w(G)
\)
holding for every nonnegative edge weighting $w$.
The forward implication follows by linearity of expectation.
Conversely, after normalizing nonzero edge weightings so that
$w(G)=1$, the finite minimax theorem yields a single probability
distribution on cuts under which every edge crosses with
probability at least $0.88672$.

In this paper, we improve the bound in
Theorem~\ref{thm:gutin-yeo-eight-elevenths} from $\frac{8}{11}$
to $\frac{15}{19}$.

\begin{restatable}{theorem}{TriangleFreeBound}
\label{thm:triangle-free-15-19}
Let $G=(V,E,w)$ be an edge-weighted triangle-free subcubic graph.
Then
\[
\operatorname{mc}(G,w)\geq\frac{15}{19}w(G).
\]
\end{restatable}

We obtain further improvements under additional odd-girth assumptions.

\begin{restatable}{theorem}{OddGirthSevenBound}
\label{thm:odd-girth-seven-47-59}
Let $G=(V,E,w)$ be an edge-weighted subcubic graph of odd girth
at least $7$. Then
\[
\operatorname{mc}(G,w)\geq\frac{47}{59}w(G).
\]
\end{restatable}

\begin{restatable}{theorem}{OddGirthNineBound}
\label{thm:odd-girth-nine-weighted}
Let $G=(V,E,w)$ be an edge-weighted subcubic graph of odd girth
at least $9$. Then
\[
\operatorname{mc}(G,w)\geq\frac45w(G).
\]
\end{restatable}

Thus Theorem~\ref{thm:odd-girth-nine-weighted} establishes
Conjecture~\ref{conj: maxcut conj} for subcubic graphs of odd girth
at least $9$. Although our bounds remain below the constant
$0.88672$ in Theorem~\ref{thm:KKV-large-girth}, they improve the
weighted guarantee under relatively small odd-girth assumptions,
allowing short even cycles.

Section~2 establishes a general lower bound
(Theorem~\ref{thm:random-left-optimize-right-triangle-free, generalized})
for weighted Max Cut in cubic graphs admitting a partition with
specified bipartite component structure.
In Section~3, we apply the general bound to an unweighted
maximum-cut partition and combine the resulting inequality with
the weight of that partition to prove
Theorem~\ref{thm:triangle-free-15-19}. Section~4 uses additional
partition structure to prove
Theorem~\ref{thm:odd-girth-seven-47-59} and
Theorem~\ref{thm:odd-girth-nine-weighted}.

\section{Setup and General Lower Bound}

Applying a strategy in~\cite{Gutin2023}, if $G$ is a subcubic graph of arbitrary odd girth, we may extend $G$ to a cubic graph with the same odd girth and weighted max-cut.

\begin{lemma}\label{lem:cubic-completion}
Let \(g\geq 3\) be odd. Every weighted subcubic graph \((G,w)\) of
odd girth at least \(g\) can be extended to a weighted cubic graph
\((\widehat G,\widehat w)\) of odd girth at least \(g\) such that
\[
w(\widehat G)=w(G)
\qquad\text{and}\qquad
\operatorname{mc}(\widehat G,\widehat w)
=\operatorname{mc}(G,w).
\]
\end{lemma}

\begin{proof}
By Erdős and Sachs~\cite{ErdosSachs1963}, bipartite cubic graphs of arbitrarily
large girth exist. Choose
a cubic bipartite graph \(B\) of girth at least \(g-1\). Subdivide an
edge \(uv\) of \(B\) with a new vertex \(x\), and denote the resulting
graph by \(B'\). Every vertex of \(B'\), except \(x\), has degree three,
while \(x\) has degree two.

For every \(s\in V(G)\), add \(3-d_G(s)\) disjoint copies of \(B'\)
and join \(s\) to the degree-two vertex \(x\) in each copy. The
resulting graph \(\widehat G\) is cubic. Each new joining edge is a
bridge, so every cycle of \(\widehat G\) lies either in \(G\) or in
one copy of \(B'\). Since \(B\) is bipartite, every odd cycle of
\(B'\) uses the subdivided edge and has length one greater than a
cycle of \(B\). It therefore has length at least \(g\). Thus
\(\widehat G\) has odd girth at least \(g\).

Finally, retain the original weights on \(E(G)\) and assign weight
zero to every new edge. Then \(w(\widehat G)=w(G)\), and cuts of
\(\widehat G\) and their restrictions to \(G\) have the same weighted
contribution. Thus 
\(
\operatorname{mc}(\widehat G,\widehat w)
=
\operatorname{mc}(G,w).
\)
\end{proof}

Throughout the remainder of this paper, let $G$ be a cubic graph with fixed odd girth. By Lemma~\ref{lem:cubic-completion}, all weighted max-cut results proved for cubic graphs extend to subcubic graphs with fixed odd girth.

Consider $(L,R)$, an unweighted max-cut partition of $G$. Since the cut is
maximum, every vertex has at least two neighbors on the opposite side, and
hence the internal edges form a matching. Write
\[
M_R:=E(G[R]),
\qquad
M_L:=E(G[L]), \qquad M = M_R \cup M_L.
\]

For $i\in\{2,3\}$, let
\[
R_i:=\{v\in R:d_{G[L,R]}(v)=i\},
\qquad
L_i:=\{v\in L:d_{G[L,R]}(v)=i\},
\qquad
V_i:=L_i\cup R_i.
\]
Thus, the indices refer to degrees in the crossing subgraph
$G[L,R]$. Partition the crossing edges as
\[
\begin{aligned}
E_{R_3L_2}&:=E(R_3,L_2),
& E_{L_3R_2}&:=E(L_3,R_2),\\
E_{22}&:=E(L_2,R_2),
& E_{33}&:=E(L_3,R_3).
\end{aligned}
\]

\begin{theorem}\label{thm:random-left-optimize-right-triangle-free, generalized}
Let $G=(V,E,w)$ be a weighted cubic graph, and let $(L,R)$ be
a partition such that $G[R]$ is a collection of bipartite components and $G[L]$ is a collection of edges and isolated vertices. Further assume that for any pair of components $C$ in $G[R]$ and $D$ in $G[L]$, $G[C \cup D]$ is bipartite. Then
\[
\operatorname{mc}(G,w)
\geq
w(G[R])
+\frac78w(G[L])
+\frac34w(E(L,R)).
\]
\end{theorem}

\begin{proof}
For each component \(C\) of
\(G[R]\) containing at least one edge, fix a bipartition
\[
V(C)=A_C\cup B_C.
\]
We will construct a new partition $(L',R')$ by randomizing the placement of $R$ and optimizing the placement of $L$. Independently for each such component, uniformly at random place one of $A_C, B_C$ in \(L'\), placing the
other class in \(R'\). For each isolated vertex of \(G[R]\), independently
choose its side in $(L',R')$ uniformly at random. Thus every edge of \(G[R]\) is crossing
and contributes $w(G[R])$.

After these random choices, place the vertices of each component of \(G[L]\)
in $(L',R')$ so as to maximize the total crossing weight of the edges incident with that
component. Since distinct components of \(G[L]\) have no edges between them,
these optimizations may be performed independently.

For an edge component \(D=uv\), we say a placement of \(D\) is
\textit{proper} if \(u\) and \(v\) lie on opposite sides of $(L',R')$,
so that the edge $uv$ is crossing.

Suppose a component \(D\) of \(G[L]\) is adjacent to \(k\leq 3\) components
\(C_1,\ldots,C_k\) of \(G[R]\). If $D$ is a singleton vertex there are two
possible placements of its vertices in \(L'\) and \(R'\). When \(D\) is an edge component, there are two proper placements of $D$ in \(L'\) and \(R'\).

Since \(G[D\cup C_i]\) is bipartite, the placement of
\(C_i\) determines a unique (proper) placement of \(D\) for which
all edges in \(E(D,C_i)\) are crossing. Call this the \textit{placement determined} by $C_i$.

\begin{claim}\label{claim:three-component-placement}
Suppose \(D\) is adjacent to \(k\leq 3\) components
\(C_1,\ldots,C_k\) of \(G[R]\). Then \(D\) can be placed so that the expected weight of the
crossing edges between \(D\) and \(C_1\cup\cdots\cup C_k\) is at least
\[
\frac34\sum_{i=1}^k w(E(D,C_i)).
\]
Further, if $D = uv$ is an edge component, the edge $uv$ is crossing.
\end{claim}

\begin{proof}[Proof of Claim~\ref{claim:three-component-placement}] Consider the three possible cases for $k$.

\begin{itemize}
    \item If \(k=1\), use the placement determined by \(C_1\). All edges (including $uv$) will be crossing.

    \item If \(k=2\), use the common (proper) placement when the placements
    determined by \(C_1\) and \(C_2\) agree. Otherwise, choose uniformly
    between the two (proper) placements of \(D\). The placement determined
    by each \(C_i\) is used with probability \(3/4\), and $uv$ will always be crossing.

    \item If \(k=3\), use the placement determined by at least two of
    the three components. For any fixed \(C_i\), its placement is not
    used only when both of the other components prefer the opposite
    placement, which occurs with probability \(1/4\).
\end{itemize}
Thus, in every case, each edge in \(E(D,C_i)\) is crossing with
probability at least \(3/4\). Linearity of expectation gives the claim.
\end{proof}

Let $D$ be the isolated vertex \(v \in L_3\). Since \(v\) has degree three, it is adjacent to at most three components of
\(G[R]\). Applying Claim~\ref{claim:three-component-placement} shows
that the expected crossing weight among the three edges incident
with \(v\) is at least three quarters of the total weight.

Next let \(D=uv\) be an edge component of \(G[L]\). Since \(G\) is cubic,
each of \(u\) and \(v\) has exactly two neighbors in \(R\). Write
\[
W_D:=w(E(D,R)).
\]

First suppose that \(D\) is adjacent to at most three components of \(G[R]\). Applying Claim~\ref{claim:three-component-placement} gives expected
external crossing weight at least
\[
\frac34 W_D.
\]
Since $uv$ is always crossing in a proper placement, the expected contribution of $D$ and its external edges is at least
\[
w(uv)+\frac34W_D
\geq
\frac78w(uv)+\frac34W_D.
\]

It remains to consider the case in which \(D = uv\) is adjacent to four distinct
components of \(G[R]\). Unlike the case when $D$ neighbors at most $3$ components, we cannot always guarantee $uv$ to be crossing. In this case, the four neighbors of \(u,v\) in \(R\)
have four independent random sides in $(L',R')$. Let \(u_1,u_2\) be the external neighbors of \(u\)
in \(R\), and let \(v_1,v_2\) be the external neighbors of \(v\) in \(R\).

We give a rule for placing \(u\) and \(v\). Call \(u\) rigid if
\(u_1,u_2\) lie on the same side in $(L',R')$ and flexible otherwise, and \(v\) analogously. If an endpoint, $u$ or $v$, is rigid, prescribe
that it be placed opposite the common side of its two external neighbors. If exactly
one endpoint is flexible, place it opposite the rigid endpoint, so that
\(uv\) crosses. If both endpoints are flexible, choose uniformly between
the two proper placements of \(u\) and \(v\).

Since the sides of \(u_1,u_2\) are independent fair choices,
\(
\Pr(u\text{ is rigid})=\Pr(u\text{ is flexible})=\frac12.
\)
When \(u\) is rigid, both edges \(uu_1,uu_2\) cross. When \(u\) is flexible,
the two neighbors lie on opposite sides, so exactly one of \(uu_1,uu_2\)
crosses. Therefore
\[
\Pr(uu_1\text{ is crossing})
=
\frac12+\frac12\cdot\frac12
=
\frac34.
\]
The same calculation applies to all four external edges. Their expected
total crossing weight is therefore
\[
\frac34W_D.
\]

For the edge \(uv\), if at least one endpoint is flexible, our rule makes
\(uv\) crossing. The rigidity events for \(u\) and \(v\) are independent,
so this occurs with probability $\frac34$.
If both endpoints are rigid, which occurs with probability \(1/4\), their
two prescribed sides are opposite with probability \(1/2\). Hence
\[
\Pr(uv\text{ is crossing})
=
\frac34+\frac14\cdot\frac12
=
\frac78.
\]
Thus in this final case the comparison rule has expected crossing weight
\[
\frac78w(uv)+\frac34W_D.
\]
Since we actually place \(u,v\) so as to maximize the weighted crossing
contribution of the component \(D\), the optimized placement performs at
least as well as the expected crossing weight. Thus, we
obtain
\[
\operatorname{mc}(G,w)
\geq
w(G[R])
+\frac78w(G[L])
+\frac34w(E(L,R)).
\]
\end{proof}
\section{Triangle Free Max Cut}
For a maximum-cut partition, both sides consist of edges and
isolated vertices, so we can apply
Theorem~\ref{thm:random-left-optimize-right-triangle-free, generalized}
in either direction.

\begin{lemma}\label{thm:random-left-optimize-right-triangle-free}
Let $G=(V,E,w)$ be a weighted triangle-free cubic graph, and let $(L,R)$ be
a max-cut partition. Then
\[
\operatorname{mc}(G,w)
\geq
\frac34 w(G)+\frac{3}{16}w(M).
\]
\end{lemma}

\begin{proof}
Since $(L,R)$ is a maximum-cut partition, every vertex has at most one
neighbor on its own side. In particular, both $G[L]$ and $G[R]$ are collections of bipartite
components. Moreover, for any component $C_L$ of $G[L]$ and any component $C_R$ of
$G[R]$, the graph $G[C_L\cup C_R]$ has at most four vertices. Since $G$ is triangle-free, $G[C_L\cup C_R]$ is bipartite. Thus Theorem~\ref{thm:random-left-optimize-right-triangle-free, generalized}
applies to $(L,R)$ and gives
\[
\operatorname{mc}(G,w)
\geq
w(G[R])
+\frac78w(G[L])
+\frac34w(E(L,R)).
\]
Applying the same theorem with the roles of $L$ and $R$ reversed gives
\[
\operatorname{mc}(G,w)
\geq
w(G[L])
+\frac78w(G[R])
+\frac34w(E(L,R)).
\]
Averaging these two inequalities, we obtain
\[
\operatorname{mc}(G,w)
\geq
\frac{15}{16}
\bigl(w(G[L])+w(G[R])\bigr)
+\frac34w(E(L,R)) =
\frac34w(G)+\frac{3}{16}w(M).
\]
\end{proof}

The original partition gives a strong bound when $w(M)$ is small,
while Lemma~\ref{thm:random-left-optimize-right-triangle-free} improves as $w(M)$ increases. Together these two bounds give an improved universal bound.

\TriangleFreeBound*

\begin{proof} By Lemma~\ref{lem:cubic-completion}, it suffices to consider
the case that \(G\) is cubic. Let $(L,R)$ be a maximum-cut partition. By Lemma~\ref{thm:random-left-optimize-right-triangle-free},
\begin{equation}\label{eq:triangle-free-random-bound}
\operatorname{mc}(G,w)
\geq
\frac34 w(G)+\frac{3}{16}w(M).
\end{equation}
On the other hand, deleting the internal edges $M$ leaves the bipartite
graph $G-M$ with bipartition $(L,R)$. Hence
\begin{equation}\label{eq:triangle-free-delete-bound}
\operatorname{mc}(G,w)\geq w(G)-w(M).
\end{equation}
Taking the convex combination of
\eqref{eq:triangle-free-random-bound} and
\eqref{eq:triangle-free-delete-bound} with coefficients
$\frac{16}{19}$ and $\frac{3}{19}$, respectively, gives
\[
\begin{aligned}
\operatorname{mc}(G,w)
&\geq
\frac{16}{19}
\left(
\frac34 w(G)+\frac{3}{16}w(M)
\right)
+
\frac{3}{19}
\left(
w(G)-w(M)
\right)\\
&=
\frac{15}{19}w(G).
\end{aligned}
\]
\end{proof}

\section{Max Cut with Higher Odd Girth}
Here we provide an improved bound for subcubic graphs with odd girth at least seven and prove Conjecture \ref{conj: maxcut conj} for graphs with odd girth at least nine. First, we prove some structural lemmas about an arbitrary max-cut $(L,R)$ of a triangle-free cubic graph.

\begin{lemma} \label{lem: degree 2 with two degree 2 nbrs}
    Let $(L, R)$ be a max-cut partition and $v \in V_2$. It is not the case that $v$ has two external neighbors $u_1, u_2 \in V_2$.
\end{lemma}

\begin{proof}
    Suppose $v$ had two external neighbors $u_1, u_2 \in V_2$. Since $G$ is triangle free, $u_1 \not\sim u_2$, thus swapping $v$ and $u_1, u_2$ makes three internal edges crossing and only loses two crossing edges. This swap increases the number of cut edges, contradicting the maximality of $(L,R)$.

\end{proof}

    \begin{figure}[htbp]
    \centering
    \includegraphics[width=0.3\linewidth]{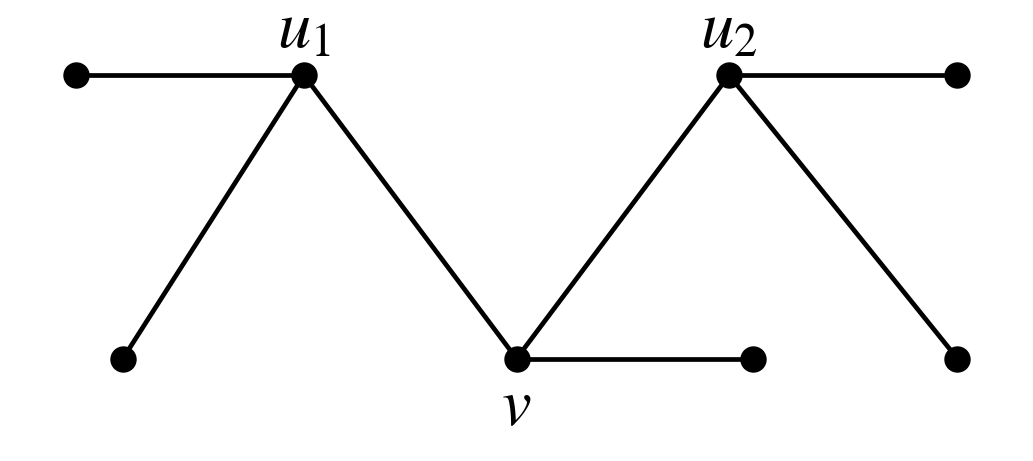}
    \caption{Configuration for Lemma~\ref{lem: degree 2 with two degree 2 nbrs}.}
    \label{fig:degree-2-two-degree-2-neighbors}
\end{figure}

\begin{lemma} \label{lem: degree 3 with two degree 2 nbrs}
    Let $(L, R)$ be a max-cut partition and $v \in V_3$. Suppose $v$ has external neighbors $u_1, u_2 \in V_2$ and $u_3 \in V_3$. Further suppose $u_1$ has an external neighbor $u_1' \in V_2$ and $u_2$ has an external neighbor $u_2' \in V_2$. Then $vu_1u_1'u_2'u_2$ forms a $5$-cycle.
\end{lemma}

\begin{proof}
    Suppose $vu_1u_1'u_2'u_2$ does not form a $5$-cycle, namely, $u_1' \not \sim u_2'$. Since $G$ is triangle-free, $u_1 \not \sim u_2$. Thus, swapping $v, u_1', u_2'$ and $u_1, u_2$ increases the number of cut edges, contradicting the maximality of $(L,R)$.
\end{proof}

\begin{figure}[htbp]
    \centering
    \includegraphics[width=0.4\linewidth]{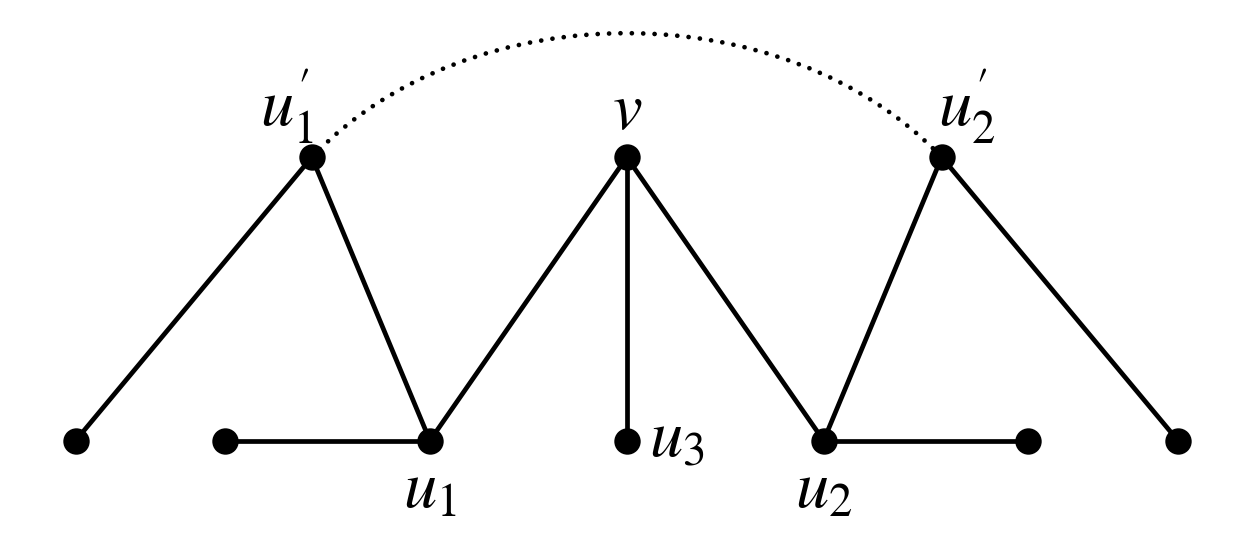}
    \caption{Configuration for Lemma~\ref{lem: degree 3 with two degree 2 nbrs}.}
    \label{fig:degree-3-two-degree-2-neighbors}
\end{figure}

\subsection{Special Max Cut Partitions} \label{subsec: special partitions}

For the rest of the paper, we consider two particular max-cut
partitions, each excluding a local
configuration on one side. Constructing the second partition from the
first allows us to bound the specific weights in terms of the first.

Begin with an arbitrary max-cut partition. Suppose that, for the
current partition $(L,R)$, there is a vertex $v\in R$ with three neighbors
$u_1,u_2,u_3\in L_2$. Swap $v$ and
$u_1,u_2,u_3$. 

\begin{figure}[htbp]
    \centering
    \includegraphics[width=0.4\linewidth]{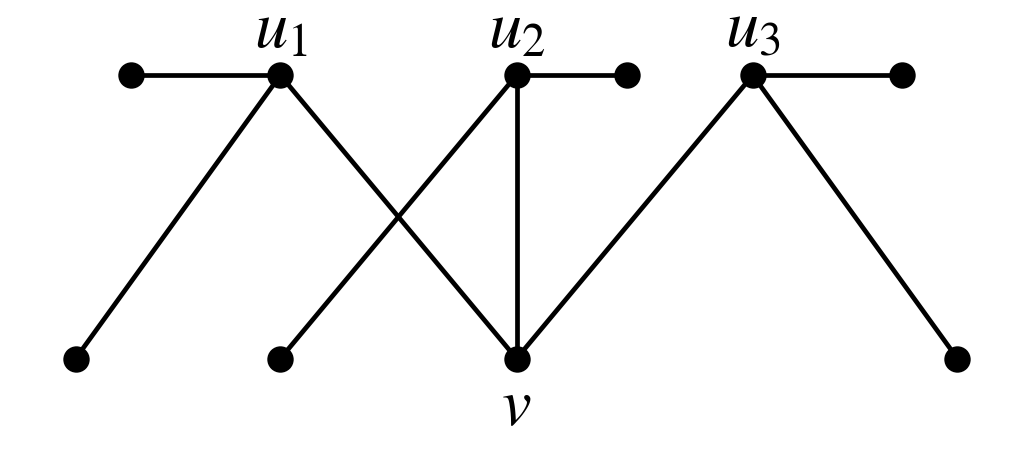}
    \caption{Swap configuration for $R$-maximality and $L$-maximality.}
    \label{fig:special-max-cut-partitions}
\end{figure}

Since $G$ is triangle-free, this operation preserves
the number of crossing edges, and hence the resulting partition is again
a maximum-cut partition. Moreover, the size of the $R$-side increases by
two. Repeat this swapping process until there are no more vertices in $R$ with three neighbors in $L_2$. If $R$ is maximal with respect to these swaps, we say $(L,R)$ is an \textit{R-maximal} max-cut partition.
Fix an $R$-maximal max-cut partition $(L,R)$. By construction, no vertex $v\in R$ has three neighbors in $L_2$.

Starting from this partition $(L,R)$, perform the analogous swapping
operation in the opposite direction. If there is a vertex $v\in L$ with
three neighbors $u_1,u_2,u_3\in R_2$, swap $v$ and
$u_1,u_2,u_3$. Again, this preserves the number of crossing edges. Repeat this swapping process until there are no more vertices in $L$ with three neighbors in $R_2$, to create a new partition $(L',R')$.  If $L'$ is maximal with respect to these swaps, we say $(L',R')$ is an \textit{L-maximal} max-cut partition.
 Fix $(L',R')$, an $L$-maximal max-cut partition resulting from the $R$-maximal max-cut partition $(L,R)$.
Then no vertex $v\in L'$ has three neighbors in $R_2'$.

During the construction of $(L',R')$, internal edges are only added to the
$L'$-side and removed from the $R'$-side. Consequently,
\begin{equation}\label{eq:MR-monotone}
M_R'\subseteq M_R \qquad\text{and}\qquad
w(M_R')\leq w(M_R).
\end{equation}

\subsection{Odd Girth at Least Seven}
 
Let $G$ be a triangle-free cubic graph with odd girth at least $7$. Consider the max-cut partitions $(L,R)$ and $(L', R')$ defined in Section~\ref{subsec: special partitions}.

Following the framework of~\cite{BerkeSzabo2007}, we will define an auxiliary conflict graph $C$ with vertex set $L_2 \subseteq L$ and a corresponding conflict graph $C'$ with vertex set $R_2' \subseteq R'$. Given $u,v \in L_2$, we say $u\sim_C v$ if $u \sim_G v$ or if $u$ and $v$ have a neighbor in the same component of $G[R]$ (which consists of $M_R$ and singleton vertices). Define $C'$ analogously for incidence of $u,v \in R_2'$ with components of $G[L']$.

An independent set in $C$ specifies vertices that can be moved
from $L$ to $R$ simultaneously so that no edge of $M_L$ has both
endpoints moved, and no component of $G[R]$ gains more than one vertex from $L$. This keeps the resulting components on the $R$-side
small, allowing us to apply
Theorem~\ref{thm:random-left-optimize-right-triangle-free, generalized}.

\begin{lemma}\label{conflict graph max degree}
The conflict graphs $C$ and $C'$ satisfy:
\begin{enumerate}
    \item $\Delta(C), \Delta(C') \leq 3$,
    \item  $C$ and $C'$ do not contain $K_4$ as a subgraph.
\end{enumerate}

This implies that $C$ and $C'$ are 3-colorable.

\end{lemma}
\begin{proof}
To prove (1), consider $v \in V(C) = L_2$. The vertex $v$ has one neighbor in $L_2$ which contributes one edge in $C$. Also, $v$ has two neighbors in $R$; let $u$ be one of them. If $u \in R_3$, then by the $R$-maximality of $(L,R)$, $u$ has at most one other neighbor in $L_2$, providing at most one more $C$-neighbor. If $u \in R_2$, then let $w$ be its unique neighbor in $R_2$. By Lemma~\ref{lem: degree 2 with two degree 2 nbrs}, $u$ has no neighbors in $L_2$ other than $v$ and $w$ has at most one neighbor in $L_2$. Thus the component of $u$ provides at most one more $C$-neighbor. Together, $d_C(v) \leq 3$.

To prove (2), suppose $C$ has a $K_4$ subgraph. Since each vertex has degree at most 3, and one neighbor comes from the matching edges of $M_L$, the $K_4$ must consist of $u_1,u_2$ and $v_1, v_2$ where $u_1u_2, v_1v_2 \in M_L$. If $u_1 \sim_C v_1, v_2$ then $u_1$ is incident with the same component as $v_1$ in $R$; call that component $C_1 \subseteq R$. Also $u_1$ is incident with the same component as $v_2$ in $R$; call that component $C_2 \subseteq R$.

If $C_1$ and $C_2$ are both singleton vertices, then $u_1\cup C_1 \cup C_2 \cup v_1v_2$ induce a $5$-cycle in $G$ which contradicts the odd girth of $G$.

If $C_1$ and $C_2$ are both edges in $M_R$, then $u_1 \in L_2$ has two neighbors in $R_2$ which contradicts Lemma~\ref{lem: degree 2 with two degree 2 nbrs}.

Thus we may assume that $C_1$ is an edge of $M_R$
and $C_2$ is an isolated vertex belonging to $R_3$.

Note that $u_2$ is incident with the same component as $v_1$ in $R$; call that component $C_1' \subseteq R$. Also $u_2$ is incident with the same component as $v_2$ in $R$; call that component $C_2' \subseteq R$. By symmetry, one of $C_1',C_2'$ consists of an edge of $M_R$,
and the other is an isolated vertex of $G[R]$ belonging to $R_3$.
To ensure $v_1$ does not violate
Lemma~\ref{lem: degree 2 with two degree 2 nbrs},
$C_1'$ must be an isolated vertex belonging to $R_3$,
and $C_2'$ must be an edge of $M_R$.

Let $C_2=\{w\}$, where $w\in R_3$ is adjacent to
$u_1,v_2\in L_2$. Since $(L,R)$ is $R$-maximal, the third
neighbor of $w$ belongs to $L_3$. Moreover, $u_1$ and $v_2$ have
neighbors $u_1',v_2'\in R_2$, respectively, in the components
$C_1$ and $C_2'$. Lemma~\ref{lem: degree 3 with two degree 2 nbrs}
therefore implies that $wu_1u_1'v_2'v_2w$ is a $5$-cycle. However, this
contradicts the assumption that $G$ has odd girth at least $7$. 

Thus $C$ contains no $K_4$.

The same argument shows that $C'$ contains no $K_4$.

By Brooks' Theorem, a $K_4$-free graph with max degree $3$ is $3$-colorable.
\end{proof}

The three-coloring lets us choose a set of vertices that can be
moved simultaneously while capturing at least one third of the
crossing weight incident with $L_2$.

\begin{lemma}
\label{lem:ML-conflict}
Let $(G,w)$ be a cubic graph with odd girth at least $7$ and let
$(L,R)$ be an $R$-maximal partition.
There exist sets $S_L \subseteq M_L$ and $S_R \subseteq E_{22} \cup E_{R_3L_2}$ such that 
\begin{equation}\label{eq:ML-delete}
\operatorname{mc}(G,w)
\geq
w(G)-w(M_L)-w(M_R)+w(S_L)-w(S_R),
\end{equation}
and
\begin{equation}\label{eq:ML-random}
\operatorname{mc}(G,w)
\geq
\frac34 w(G)
+\frac18w(M_L)
+\frac14w(M_R)
-\frac18w(S_L)
+\frac14w(S_R).
\end{equation}
Moreover, the sets $S_L$ and $S_R$ may be chosen so that
\begin{equation}\label{eq:ML-SR}
w(S_R)
\geq
\frac13
\left(
w(E_{22})+w(E_{R_3L_2})
\right).
\end{equation}
\end{lemma}

\begin{proof}
Consider $C$, the conflict graph for $L_2$. By Lemma~\ref{conflict graph max degree}, the conflict graph \(C\) is
3-colorable. Let $L_2=X_1\cup X_2\cup X_3$ be a partition into three independent sets of \(C\).

Since $X_1,X_2,X_3$ partition $L_2$ and every edge of
$E_{22}\cup E_{R_3L_2}$ has exactly one endpoint in $L_2$,
some color class $X_i$ is incident with at least one third
of the total weight of these edges.

Let $X:=X_i$, and let $S_R$ consist of the edges of
$E_{22}\cup E_{R_3L_2}$ incident with $X$. Then
\[
w(S_R)\geq
\frac13\left(w(E_{22})+w(E_{R_3L_2})\right),
\]
which proves~\eqref{eq:ML-SR}.
Also let $S_L$ be the edges of $M_L$ incident with $X$.

Now move every vertex of \(X\) from \(L\) to \(R\), and denote the
resulting partition by \((L^*,R^*)\). Since \(X\) is an independent set in the
conflict graph $C$, no edge of \(M_L\) has both endpoints in \(X\). Hence the
edges of \(S_L\) become crossing, while the remaining edges of \(M_L\)
remain internal to \(L^*\). Thus
\[
E(G[L^*])=M_L\setminus S_L.
\]
Similarly, the crossing edges which become internal to
\(R^*\) are the edges of \(S_R\), so
\[
E(G[R^*])=M_R\cup S_R.
\]
The cut \((L^*,R^*)\) contains every
edge of \(G\) except the edges of
\(
(M_L\setminus S_L)\cup M_R\cup S_R.
\)
Therefore
\[
\begin{aligned}
\operatorname{mc}(G,w)
&\geq
w(G)-w(M_L\setminus S_L)-w(M_R)-w(S_R)\\
&=
w(G)-w(M_L)-w(M_R)+w(S_L)-w(S_R),
\end{aligned}
\]
which is~\eqref{eq:ML-delete}.

It remains to prove~\eqref{eq:ML-random}. We claim that the partition
\((L^*,R^*)\) satisfies the hypotheses of
Theorem~\ref{thm:random-left-optimize-right-triangle-free, generalized}. First, \(G[L^*]\) is a collection of edges and isolated vertices, since
\(M_L\setminus S_L\) is a matching.

Next consider \(G[R^*]\). Before moving \(X\), every component of \(G[R]\)
was either an edge of \(M_R\) or an isolated vertex. Since \(X\) is
independent in the conflict graph, no component of \(G[R]\) is adjacent
to two distinct vertices of \(X\). Hence every component of \(G[R^*]\)
contains at most one vertex of \(X\).

Moreover, a vertex \(x\in X\subseteq L_2\) cannot have both of its
neighbors in \(R_2\), by
Lemma~\ref{lem: degree 2 with two degree 2 nbrs}. Among the
two components of \(G[R]\) adjacent to \(x\), at most one is an edge of
\(M_R\).
Therefore every component of \(G[R^*]\) has order at most $4$. Since \(G\) has odd girth at least \(7\), every such component is
bipartite.

Finally, let \(C_R\) be a component of \(G[R^*]\) and \(C_L\) a component
of \(G[L^*]\). We have $|V(C_R\cup C_L)|\leq 6$, which implies $G[C_R \cup C_L]$ is bipartite since $G$ has odd girth at least $7$.

Applying Theorem~\ref{thm:random-left-optimize-right-triangle-free, generalized} to \((L^*,R^*)\) gives
\[
\begin{aligned}
\operatorname{mc}(G,w)
&\geq
w(G[R^*])
+\frac78w(G[L^*])
+\frac34w(E(L^*,R^*))\\
&=
w(M_R)+w(S_R)
+\frac78\bigl(w(M_L)-w(S_L)\bigr)\\
&\qquad
+\frac34\Bigl(
w(G)-w(M_R)-w(S_R)-w(M_L)+w(S_L)
\Bigr).
\end{aligned}
\]
\[
\operatorname{mc}(G,w)
\geq
\frac34 w(G)
+\frac18w(M_L)
+\frac14w(M_R)
-\frac18w(S_L)
+\frac14w(S_R),
\]
which is~\eqref{eq:ML-random}.
\end{proof}

We now apply the same construction to $(L',R')$, with the roles
of the two sides reversed. This time we choose the color class
so that the edges of $M_R'$ that become crossing have large total
weight. Each edge of $M_R'$ meets two of the three color classes,
so some class captures at least $\frac23w(M_R')$.

\begin{lemma}
\label{lem:MRprime-conflict}
Let $(G,w)$ be a cubic graph of odd girth at least $7$, and let
$(L',R')$ be the $L$-maximal partition derived from $(L,R)$. Then there exist sets
\[
S_R' \subseteq M_R'
\qquad\text{and}\qquad
S_L' \subseteq E_{22}'\cup E_{L_3R_2}'
\]
such that
\begin{equation}\label{eq:MRprime-delete}
\operatorname{mc}(G,w)
\geq
w(G)-w(M_L')-w(M_R')+w(S_R')-w(S_L'),
\end{equation}
and
\begin{equation}\label{eq:MRprime-random}
\operatorname{mc}(G,w)
\geq
\frac34 w(G)
+\frac14w(M_L')
+\frac18w(M_R')
+\frac14w(S_L')
-\frac18w(S_R').
\end{equation}
Moreover, the sets $S_L'$ and $S_R'$ may be chosen so that
\begin{equation}\label{eq:MRprime-SR}
w(S_R')
\geq
\frac23 w(M_R').
\end{equation}
\end{lemma}

\begin{proof}
Consider the conflict graph $C'$ for $R_2'$. By
Lemma~\ref{conflict graph max degree}, it is $3$-colorable; let
\(
R_2'=X_1\cup X_2\cup X_3
\)
be a partition into three independent sets.

Since the endpoints of every edge of $M_R'$ are adjacent in $C'$,
they belong to distinct color classes. Thus every edge of $M_R'$
is incident with exactly two of the three color classes, so some
color class $X_i$ is incident with edges of $M_R'$ of total weight
at least $\frac23w(M_R')$.

Set $X:=X_i$, and let $S_R'$ be the set of edges of $M_R'$
incident with $X$. Then
\[
w(S_R')\geq \frac23w(M_R').
\]
Let $S_L'$ be the set of edges of
$E_{22}'\cup E_{L_3R_2}'$ incident with $X$.

Move every vertex of $X$ from $R'$ to $L'$, and denote the resulting
partition by $(L^*,R^*)$. Then
\[
E(G[R^*])=M_R'\setminus S_R'
\quad\text{and}\quad
E(G[L^*])=M_L'\cup S_L'.
\]
The weight of the cut $(L^*,R^*)$ gives
\[
\operatorname{mc}(G,w)
\geq
w(G)-w(M_L')-w(M_R')+w(S_R')-w(S_L'),
\]
proving~\eqref{eq:MRprime-delete}.

As in the proof of Lemma~\ref{lem:ML-conflict}, the independence of $X$
in $C'$ ensures that
Theorem~\ref{thm:random-left-optimize-right-triangle-free, generalized}
applies, with the roles of $L^*$ and $R^*$ reversed. Therefore
\[
\begin{aligned}
\operatorname{mc}(G,w)
&\geq
w(G[L^*])+\frac78w(G[R^*])
+\frac34w(E(L^*,R^*))\\
&=
\frac34 w(G)
+\frac14w(M_L')
+\frac18w(M_R')
+\frac14w(S_L')
-\frac18w(S_R'),
\end{aligned}
\]
proving~\eqref{eq:MRprime-random}.
\end{proof}

The next bound exploits a different tradeoff; we can make every
edge of $M_L$ crossing while retaining at least half the
crossing weight incident with $L_2$. This is useful when substantial weight lies
in $M_L$.

\begin{lemma}
\label{lem:half-weight}
Let $(G,w)$ be a cubic graph, and let $(L,R)$
be a maximum-cut partition of $G$.
Then
\begin{equation}\label{eq:half-weight}
\operatorname{mc}(G,w)
\geq
w(M_L)+w(E_{33})+w(E_{L_3R_2})
+\frac12w(E_{22})
+\frac12w(E_{R_3L_2}).
\end{equation}
\end{lemma}

\begin{proof}
Starting from $(L,R)$, independently for each edge $uv\in M_L$, choose one
of $u,v$ uniformly at random and move the chosen endpoint from $L$ to $R$.
Let $(L^*,R^*)$ denote the resulting random partition.

Exactly one endpoint of every edge of $M_L$ is moved, so every edge of
$M_L$ is crossing in $(L^*,R^*)$. Thus these edges contribute
\(
w(M_L).
\)

No endpoint of an edge in $E_{33}$ or $E_{L_3R_2}$ belongs to $L_2$.
Therefore none of their endpoints is moved, and every such edge remains
crossing. These classes contribute
\(
w(E_{33})+w(E_{L_3R_2}).
\)

Now let $e\in E_{22}\cup E_{R_3L_2}$. The edge $e$ remains
crossing precisely when its $L_2$ endpoint is not chosen to move.
The expected crossing weight contributed by these two edge
classes is
\(
\frac12w(E_{22})+\frac12w(E_{R_3L_2}).
\)

Hence some realization of $(L^*,R^*)$ has crossing weight at least this
expectation. Therefore
\[
\operatorname{mc}(G,w)
\geq
w(M_L)+w(E_{33})+w(E_{L_3R_2})
+\frac12w(E_{22})
+\frac12w(E_{R_3L_2}).
\]
\end{proof}

Now we combine the odd girth at least 7 bounds to get a universal bound for $\operatorname{mc}(G,w)$ in the odd girth at least 7 setting.

\OddGirthSevenBound*

\begin{proof}By Lemma~\ref{lem:cubic-completion}, it suffices to consider
the case that \(G\) is cubic. Let $(L,R)$ be an $R$-maximal maximum-cut partition, and let
$(L',R')$ be an $L$-maximal partition derived from $(L,R)$.

From Lemma~\ref{lem:ML-conflict}, taking $3$ times
\eqref{eq:ML-delete} plus $24$ times \eqref{eq:ML-random} gives
\[
\begin{aligned}
27\operatorname{mc}(G,w)
&\geq 21w(G)+3w(M_R)+3w(S_R).
\end{aligned}
\]
Then applying \eqref{eq:ML-SR} gives
\begin{equation}\label{eq:odd7-first-reduced}
\begin{aligned}
27\operatorname{mc}(G,w)
&\geq 21w(G)+3w(M_R)
+w(E_{22})+w(E_{R_3L_2}).
\end{aligned}
\end{equation}

From Lemma~\ref{lem:MRprime-conflict}, taking $6$ times
\eqref{eq:MRprime-delete} plus $24$ times
\eqref{eq:MRprime-random} gives
\[
\begin{aligned}
30\operatorname{mc}(G,w)
&\geq 24w(G)-3w(M_R')+3w(S_R').
\end{aligned}
\]
Then applying \eqref{eq:MRprime-SR} followed by
\eqref{eq:MR-monotone} gives
\begin{equation}\label{eq:odd7-second-reduced}
\begin{aligned}
30\operatorname{mc}(G,w)
&\geq 24w(G)-w(M_R')\\
&\geq 24w(G)-w(M_R).
\end{aligned}
\end{equation}

Finally, since
$M_L,M_R,E_{22},E_{33},E_{L_3R_2},E_{R_3L_2}$
partition $E(G)$, $2$ times \eqref{eq:half-weight} from
Lemma~\ref{lem:half-weight} can be written as
\begin{equation}\label{eq:odd7-half-weight-doubled}
2\operatorname{mc}(G,w)
\geq
2w(G)-2w(M_R)-w(E_{22})-w(E_{R_3L_2}).
\end{equation}
Adding \eqref{eq:odd7-first-reduced},
\eqref{eq:odd7-second-reduced}, and
\eqref{eq:odd7-half-weight-doubled}, all other terms cancel, yielding
\[
59\operatorname{mc}(G,w)\geq 47w(G).
\]
\end{proof}

\subsection{Odd Girth at Least Nine}

A structural lemma of Berke and Szabó gives a natural partition to apply Theorem~\ref{thm:random-left-optimize-right-triangle-free, generalized} for odd girth at least nine.

\begin{lemma}[Berke and Szabó~\cite{BerkeSzabo2007}]
\label{lem:berke-szabo}
Let \(G\) be a triangle-free subcubic graph.
Then there exists a partition
\[
V(G)=I\cup B
\]
such that \(I\) is an independent set and every component of \(G[B]\) has
order at most \(6\).
\end{lemma}

\OddGirthNineBound*

\begin{proof}
By Lemma~\ref{lem:cubic-completion}, it suffices to consider
the case that \(G\) is cubic. By Lemma~\ref{lem:berke-szabo}, there exists a partition
\(
V(G)=I\cup B
\)
such that \(I\) is an independent set and every component of \(G[B]\) has
order at most \(6\).

The partition \((I,B)\) satisfies the hypotheses of
Theorem~\ref{thm:random-left-optimize-right-triangle-free, generalized}. First, \(G[I]\) consists entirely of isolated
vertices. Moreover, every component \(C\) of \(G[B]\) has at most $6$ vertices. Since \(G\) has odd girth at
least \(9\), $C$ is bipartite. Finally, let \(C\) be a component of \(G[B]\) and let \(v\in I\). Since
\(
|V(C\cup\{v\})|\leq 7,
\) odd girth implies
\(G[C\cup\{v\}]\) is bipartite.

By
Theorem~\ref{thm:random-left-optimize-right-triangle-free, generalized},
\begin{equation}\label{eq:odd9-generalized-bound}
\operatorname{mc}(G,w)
\geq
w(G[B])+\frac34 w(E(I,B)).
\end{equation}

On the other hand, the partition \((I,B)\) itself cuts every edge of
\(E(I,B)\), so
\begin{equation}\label{eq:odd9-partition-bound}
\operatorname{mc}(G,w)\geq w(E(I,B)).
\end{equation}

Taking \(\frac15\) times \eqref{eq:odd9-partition-bound} and
\(\frac45\) times \eqref{eq:odd9-generalized-bound} gives
\[
\begin{aligned}
\operatorname{mc}(G,w)
&\geq
\frac15 w(E(I,B))
+\frac45\left(
w(G[B])+\frac34w(E(I,B))
\right)\\
&=
\frac45 w(G[B])
+\frac45 w(E(I,B))\\
&=
\frac45 w(G).
\end{aligned}
\]
\end{proof}

\section{Conclusion}

Define
\(
b(G)=\frac{\operatorname{mc}(G)}{|E(G)|} \) and \(b(G,w)=\frac{\operatorname{mc}(G,w)}{w(G)}.
\)
It is natural to ask whether
\(b(G,w)\geq b(G)\) for every triangle-free subcubic graph \(G\)
and every nonnegative edge weighting \(w\). This is false even for
cubic graphs. Consider the Wagner graph \(W_8\), obtained from an \(8\)-cycle by
adding the four edges joining opposite vertices. This graph is
triangle-free and cubic, and its maximum cut has size \(10\), so
\(
b(W_8)=\frac{10}{12}=\frac56.
\)
However, \(W_8\) contains a \(5\)-cycle \(C\). Assign weight \(1/5\)
to each edge of \(C\) and weight zero to every other edge.
Every cut misses at least one edge of \(C\), and a cut crossing
four of its edges extends to a cut of \(W_8\). Hence
\(
b(W_8,w)=\frac45<\frac56=b(W_8).
\)

The Wagner graph example shows that weights can lower the maximum-cut ratio
of an individual graph. However, the resulting ratio is still
$4/5$, the unweighted bound for the entire triangle-free
subcubic class. Thus the example leaves open whether allowing
weights changes the worst-case ratio over an entire odd-girth
class.

For each odd integer \(g\geq5\), let \(\mathcal{G}_g\) be the class
of finite simple subcubic graphs with at least one edge and odd
girth at least \(g\), including bipartite graphs. Define
\[
b_g^{\mathrm{unw}}
=\inf_{G\in\mathcal{G}_g} b(G) \qquad \text{and} \qquad
b_g^{\mathrm{wt}}
=\inf_{G\in\mathcal{G}_g}
  \inf_{\substack{w:E(G)\to[0,\infty)\\w(G)>0}} b(G,w).
\]
Uniform edge weights show that
\(b_g^{\mathrm{wt}}\leq b_g^{\mathrm{unw}}\).

\begin{question}\label{q:weighted-unweighted}
For every odd integer \(g\geq5\), is
\(
b_g^{\mathrm{wt}}=b_g^{\mathrm{unw}}?
\)
\end{question}

Our results give weighted lower bounds under stronger odd-girth
assumptions, but determining the sharp bounds remains open even
in the unweighted setting. As discussed in the introduction,
Zhu~\cite{Zhu2009} proved that every $2$-connected triangle-free
subcubic graph outside the seven exceptional graphs has
maximum-cut ratio at least $17/21$. Since each exceptional graph
contains a $5$-cycle and maximum-cut sizes can only increase over bridges,
this bound holds for every subcubic graph of odd girth at
least seven. Together with
Theorems~\ref{thm:odd-girth-seven-47-59}
and~\ref{thm:odd-girth-nine-weighted},
\[
\begin{aligned}
\frac{17}{21}
&\leq b_7^{\mathrm{unw}},
&
\frac{47}{59}
&\leq b_7^{\mathrm{wt}},\\
\frac{17}{21}
&\leq b_9^{\mathrm{unw}},
&
\frac45
&\leq b_9^{\mathrm{wt}}.
\end{aligned}
\]

\begin{question}\label{q:sharp-unweighted}
For each odd integer \(g\geq7\), what is
\(b_g^{\mathrm{unw}}\)?
In particular, what are \(b_7^{\mathrm{unw}}\) and
\(b_9^{\mathrm{unw}}\)?
\end{question}

Improving any of these lower bounds would be interesting.

Finally, Theorem~\ref{thm:random-left-optimize-right-triangle-free, generalized}
suggests a structural approach to the weighted \(4/5\) conjecture.

\begin{question}\label{q:bipartite-components}
Does every triangle-free subcubic graph \(G\) admit a partition
\(V(G)=L\cup R\) such that \(L\) is independent and, for every
component \(C\) of \(G[R]\) and every vertex \(v\in L\), the
graph \(G[V(C)\cup\{v\}]\) is bipartite?
\end{question}

One possible approach is to strengthen the component-order bound
in Lemma~\ref{lem:berke-szabo}, taking the odd girth into account.

\begin{question}\label{q:small-components-five}
Does every subcubic graph \(G\) of odd girth exactly five admit
a partition \(V(G)=L\cup R\) such that \(L\) is independent and
every component of \(G[R]\) has order at most three?
\end{question}

\begin{question}\label{q:small-components-seven}
Does every subcubic graph \(G\) of odd girth at least seven admit
a partition \(V(G)=L\cup R\) such that \(L\) is independent and
every component of \(G[R]\) has order at most five?
\end{question}

An affirmative answer to Question~\ref{q:bipartite-components}
would imply the weighted \(4/5\) conjecture. Likewise, an
affirmative answer to Question~\ref{q:small-components-five}
or Question~\ref{q:small-components-seven} would yield the weighted \(4/5\) conjecture
for the corresponding odd-girth class.

\section{Use of AI}
The main mathematical ideas were developed solely by the author, including the structural analysis of maximum-cut partitions, the randomized cut construction underlying Theorem~\ref{thm:random-left-optimize-right-triangle-free, generalized}, and the use of odd-girth restrictions to obtain stronger bounds.
ChatGPT-6 was used to simplify the component inequalities and compute the convex combinations in the proof of Theorem~1.5.
It also assisted with improving LaTeX formatting and checking the manuscript for potential mathematical errors.
The author reviewed and verified all AI-assisted material and takes full responsibility for the content of this paper.

\bibliographystyle{plain}  
\bibliography{references} 

@article{Hopkins1982,
  author    = {G. Hopkins and W. Staton},
  title     = {Extremal Bipartite Subgraphs of Cubic Triangle-Free Graphs},
  journal   = {Journal of Graph Theory},
  volume    = {6},
  number    = {2},
  pages     = {115--121},
  year      = {1982},
  doi       = {10.1002/jgt.3190060113}
}

@article{XuYu2008,
  author  = {Xu, Baogang and Yu, Xingxing},
  title   = {Triangle-free subcubic graphs with minimum bipartite density},
  journal = {Journal of Combinatorial Theory, Series B},
  volume  = {98},
  number  = {3},
  pages   = {516--537},
  year    = {2008},
  issn    = {0095-8956},
  doi     = {10.1016/j.jctb.2007.09.001},
  url     = {https://www.sciencedirect.com/science/article/pii/S0095895607001062}
}

@article{Zhu2009,
  author  = {Zhu, Xuding},
  title   = {Bipartite subgraphs of triangle-free subcubic graphs},
  journal = {Journal of Combinatorial Theory, Series B},
  volume  = {99},
  number  = {1},
  pages   = {62--83},
  year    = {2009},
  doi     = {10.1016/j.jctb.2008.04.005}
}

@article{Yannakakis1978,
  author    = {M. Yannakakis},
  title     = {Node- and Edge-Deletion {NP}-Complete Problems},
  journal   = {Proceedings of the Tenth Annual ACM Symposium on Theory of Computing},
address = {New York, NY, USA},
  volume    = {1},
  number    = {1},
  pages     = {253--264},
  year      = {1978},
  doi       = {10.1145/800133.804355}
}

@article{PoljakTurzik1986,
  author  = {Poljak, Svatopluk and Turz{\'i}k, Daniel},
  title   = {A polynomial time heuristic for certain subgraph optimization problems with guaranteed worst case bound},
  journal = {Discrete Mathematics},
  volume  = {58},
  number  = {1},
  pages   = {99--104},
  year    = {1986},
  issn    = {0012-365X},
  doi     = {10.1016/0012-365X(86)90192-5},
  url     = {https://www.sciencedirect.com/science/article/pii/0012365X86901925}
}

@article{HalperinLivnatZwick2004,
  author  = {Halperin, Eran and Livnat, Dror and Zwick, Uri},
  title   = {{MAX CUT} in cubic graphs},
  journal = {Journal of Algorithms},
  volume  = {53},
  number  = {2},
  pages   = {169--185},
  year    = {2004},
  issn    = {0196-6774},
  doi     = {10.1016/j.jalgor.2004.06.001},
  url     = {https://www.sciencedirect.com/science/article/pii/S0196677404000902}
}

@article{KardosKralVolec2012,
  author  = {Kardo{\v{s}}, Franti{\v{s}}ek and Kr{\'a}{\v{l}}, Daniel and Volec, Jan},
  title   = {Maximum edge-cuts in cubic graphs with large girth and in random cubic graphs},
  journal = {Random Structures \& Algorithms},
  volume  = {41},
  number  = {4},
  pages   = {506--520},
  year    = {2012},
  issn    = {1042-9832},
  doi     = {10.1002/rsa.20471},
  url     = {https://doi.org/10.1002/rsa.20471}
}

@article{BondyLocke1986,
  author  = {Bondy, J. A. and Locke, S. C.},
  title   = {Largest bipartite subgraphs in triangle-free graphs with maximum degree three},
  journal = {Journal of Graph Theory},
  volume  = {10},
  number  = {4},
  pages   = {477--504},
  year    = {1986},
  issn    = {0364-9024},
  doi     = {10.1002/jgt.3190100407},
  url     = {https://doi.org/10.1002/jgt.3190100407}
}

@article{Gutin2023,
  author    = {G. Gutin and A. Yeo},
  title     = {Lower Bounds for Maximum Weighted Cut},
  journal   = {SIAM Journal on Discrete Mathematics},
  volume    = {37},
  number    = {2},
  pages     = {1142--1161},
  year      = {2023},
  doi       = {10.1137/21M1411913}
}

@article{BerkeSzabo2007,
  author  = {Robert Berke and Tibor Szab{\'o}},
  title   = {Relaxed two-coloring of cubic graphs},
  journal = {Journal of Combinatorial Theory, Series B},
  volume  = {97},
  number  = {4},
  pages   = {652--668},
  year    = {2007},
  doi     = {10.1016/j.jctb.2006.12.001}
}

@inproceedings{BermanKarpinski1999,
  author    = {Piotr Berman and Marek Karpinski},
  title     = {On Some Tighter Inapproximability Results
               (Extended Abstract)},
  booktitle = {Automata, Languages and Programming},
  series    = {Lecture Notes in Computer Science},
  volume    = {1644},
  pages     = {200--209},
  publisher = {Springer},
  year      = {1999},
  doi       = {10.1007/3-540-48523-6_17}
}

@article{ErdosSachs1963,
  author  = {Paul Erd{\H{o}}s and Horst Sachs},
  title   = {Regul{\"a}re Graphen gegebener Taillenweite mit minimaler
             Knotenzahl},
  journal = {Wiss. Z. Martin-Luther-Univ. Halle-Wittenberg
             Math.-Natur. Reihe},
  volume  = {12},
  pages   = {251--257},
  year    = {1963}
}
\end{document}